\documentstyle[12pt]{article}

  \newcommand{\const}{\rm const}
  
  \newcommand{\supp}{\rm supp}

  \newcommand{\vraisup}{\rm vraisup}
  \newcommand{\Dom}{\rm  Dom}

\begin{document}

   \begin{center}

    {\bf A note about associate and dual spaces }\\

  \vspace{4mm}

  {\bf  to the Grand Lebesgue ones } \\

\vspace{7mm}

  {\bf   Ostrovsky E., Sirota L.}\\

\vspace{4mm}

 Israel,  Bar-Ilan University, department of Mathematic and Statistics, 59200, \\

\vspace{4mm}

e-mails: \  eugostrovsky@list.ru \\
sirota3@bezeqint.net \\

\vspace{4mm}

  {\bf Abstract} \\

\vspace{4mm}

 \end{center}

 \ \  We  describe in this short  article the associate and dual (conjugate) spaces to the Grand Lebesgue Spaces
 by means of its embedding to the suitable exponential Orlicz ones. \par

\vspace{4mm}


\vspace{4mm}


\vspace{4mm}

 \section{ Definitions. Notations. Previous results. Statement of problem.}

 \vspace{4mm}

\begin{center}

 \ {\it A.  Grand Lebesgue \ $ \  G(\psi)  \  $ spaces.} \par

 \end{center}

\vspace{4mm}

 \  Let $ \  (X,B, \mu)  \  $ be non-trivial measurable space with separable diffuse measure  $ \ \mu. \ $ This imply that
 for arbitrary set $ \ A \in B, \ $  for which $ \  \mu(A) \in (0,\infty) \ $ there exists a subset $ \ A_1 \subset A \ $
 such that $ \ \mu(A_1) = \mu(A)/2. \ $ \par
 \ The {\it separability} is understood relative the usually distance function

$$
\rho(A_1,A_2) \stackrel{def}{=} \mu(A_1 \setminus A_2) + \mu(A_2 \setminus A_1), \ A_1, A_2 \in B.
$$

\ Note in particular that despite these measures are atomless, but the random values (measurable functions) defined on these spaces may have a discrete distribution. \par

  \ Let also now  $  \psi = \psi(p), \ p \in (a,b), \ a = \const \in [1,\infty), \ b = \const \in (a, \infty] $ be certain bounded
from below: $ \inf \psi(p) > 0 $ continuous inside the semi - open interval  $ p \in (a,b) $
numerical valued function. We can and will suppose  without loss of generality

 $$
  \inf_{p \in (a,b)} \psi(p) = 1 \eqno(1.0)
$$
and $ \ a = \inf \{ p, \ \psi(p) < \infty \}; \ \hspace{4mm}  \ b = \sup \{ p, \ \psi(p) < \infty  \}, $  so that
$ \supp \ \psi = [a, b) $ or  $ \supp \ \psi = [a, b] \ $  or
  $ \supp \ \psi = (a, b) $  or in turn $ \supp \ \psi = (a, b]. $ The set of all such a functions will be
denoted by  $ \ \Psi(a,b) = \{ \psi(\cdot)  \}; \ $

$$
\ \Psi := \cup_{(a,b):  \ 1 \le a < b \le \infty}  \Psi(a,b).
$$

 \ One can define formally $ \ \psi(p) = +\infty, \ p \notin \supp \ \psi \ $ and $ \  C/\infty := 0. \ $ \par

 \ Of course, in the case when the measure  $ \ \mu(X) \ $ is finite; one can accept without loss of generality  $ \ \mu(X) = 1, \ $
 (a probabilistic case), then it is reasonable to admit only $ \ \supp \ \psi = [1,b) \ $ or $ \ \supp \ \psi = [1,b]; \ $ the last restriction may
 holds true iff $ \ b < \infty. \ $ \par

\vspace{4mm}

 {\bf Definition 1.1.} (See [4], [5]-[7], [11]-[12], [13], [10].) \par

 \ By definition, the (Banach) Grand Lebesgue Space (GLS)    $  \ G \psi  = G\psi(a,b) $
consists on all the real (or complex) numerical valued measurable functions
(random variables, r.v., in the case when $ \ \mu(X) = 1)   \hspace{5mm}   \  f: \ X \to R \ $  defined on our measurable space and having a finite norm

$$
|| \ f \ || = ||f||G\psi \stackrel{def}{=} \sup_{p \in (a,b)} \left[ \frac{|f|_p}{\psi(p)} \right]. \eqno(1.1)
$$

 \ Here and in what follows the notation $ \  |f|_p  \ $ denotes an ordinary Lebesgue-Riesz $ \ L_p = L_p(X,\mu) \ $ norm for the measurable function $ \ f: \ $

 $$
 |f|_p  \stackrel{def}{=} \left[ \ \int_X |f(x)|^p \ \mu(dx) \ \right]^{1/p}, \ p \ge 1.
 $$

\vspace{4mm}

 \ The function $ \  \psi = \psi(p) \  $ is said to be  the {\it  generating function } for this space. \par

 \ If for instance $ \ \mu(X) = 1 \ $ and $ \ \psi(p) = \psi_{r}(p) = 1, \ p \in [1,r],  \  $ where $ \ r = \const \in [1,\infty), $ (an extremal case), then the correspondent
 $ \  G\psi_{r}(p)  \  $ space coincides  with the classical Lebesgue-Riesz space $ \ L_r(X,\mu) = L_r(X) = L_r \ $

$$
||\xi|| G\psi_{r} = |\xi|_r, \ r \in [1, \infty).
$$

\vspace{4mm}

\  Furthermore,  let now $  \eta = \eta(z), \ z \in S $ be arbitrary family of random variables defined on any set $ \ z \in S \ $ such that

$$
\exists a, b = \const\in (1,\infty], \ a < b, \ \forall p \in [1,b)  \ \Rightarrow  \psi_S(p) := \sup_{z \in S} |\eta(z)|_p  < \infty.
$$
 \ The function $  p \to \psi_S(p)  $ is named as a {\it  natural} function for the  family  of random variables $  S.  $  Obviously,

$$
\sup_{z \in S} ||\eta(z)||G\Psi_S = 1.
$$

 \ The family $ \ S \ $ may consists on the unique r.v., say $  \  \Delta: \ $

$$
\psi_{\Delta}(p):= |\Delta|_p,
$$
if of course  the last function is finite for some values $ \  p \in( p_1,p_2); \ p_0  > 1, \ p_1 \in (p_0, \infty]. \  $\par

 \ These spaces appearers at first in [15]; see also [11]-[13], [4], [6]-[7], [10], [20] etc.
  They are applied in statistics, theory of random processes and fields, theory of Partial Differential
 Equations, theory of operators and so one. \par

\vspace{4mm}

\begin{center}

 \ {\it B. \ $ \ B(\phi) \ $ spaces. } \par

\end{center}

 \ A very important subclass of these spaces form the so - called $ \  B(\phi) \ $ spaces.
Let now  $  \ (\Omega = \{\omega\}, F, {\bf P} ) \ $  be certain sufficiently rich probability space.
Let also $ \ \phi = \phi(\lambda), \ \lambda \in (-\lambda_0, \ \lambda_0), \ \lambda_0 = \const \in (0, \infty] $
be certain even strong convex  which takes positive values for positive arguments twice continuous differentiable function, briefly: Young-Orlicz function,
such that

$$
\phi(0) = 0; \ \phi^{''}(0) \in (0,\infty).
$$

 \ For instance: $ \ \phi(\lambda) = \phi_2(\lambda) = \lambda^2/2;  \ $ is the so-called subgaussian case.\par

 \ We denote the set of all these Young-Orlicz function as $ \ \Phi  = \Phi_{\lambda_0}= \{ \ \phi \ \}. \ $\par

\vspace{4mm}

{\bf Definition 1.2.} (See [15], [11]-[12], [13].)\par

 \ We will say by definition that the centered numerical valued random variable (r.v) $  \ \xi \ $ belongs
to the space $  B(\phi), \ \phi \in \Phi,\ $ if there exists certain non-negative constant $ \ \tau \ge 0 \ $ such that

$$
\forall \lambda \in (- \lambda_0, \ \lambda_0) \ \Rightarrow  {\bf E}\exp(\lambda \ \xi) \le \exp(\phi(\lambda \ \tau)). \eqno(1.2).
$$

 \ The minimal non-negative value $ \ \tau \ $ satisfying (1.2) for all the values $ \  \lambda \in (- \lambda_0, \ \lambda_0), \ $
 is named a $ \  B(\phi) \ $ norm of the variable $ \ \xi, \ $ write

$$
||\xi||B(\phi) \stackrel{def}{=} \inf \{\tau, \tau > 0: \  \forall \lambda \in (- \lambda_0, \ \lambda_0) \ \Rightarrow  {\bf E}\exp(\lambda \ \xi) \le \exp(\phi(\lambda \ \tau)) \}. \eqno(1.3)
$$

\vspace{4mm}

 \ These spaces are very convenient for the investigation of the r.v. having a exponential decreasing tail of distribution, for instance, for investigation of the limit
theorem, the exponential bounds of distribution for sums of random variables, non-asymptotical properties, problem of continuous and weak compactness of random
fields, study of Central Limit Theorem in the Banach space etc. The detail investigation of these spaces may be found in [15], [4], [6]-[7], [8], [13]-[14], [16], [17], [20]. \par

\ One unexpected new application of these spaces  may be found in a recent  article [1]. \par

 \ The space $ \  B(\phi), \ \phi \in \Phi \ $ with respect to the norm  $ || \xi ||B(\phi) \ $ and ordinary algebraic operations is a rearrangement invariant Banach space in the classical sense
  [3], chapters 1,2; which is in turn isomorphic to the subspace of the space $ \ G\psi_{\phi}, \ $ where

 $$
 \psi_{\phi}(p) := \frac{p}{\phi^{-1}(p)}.
 $$
  consisting on all the centered  r.v. having a finite norm $ \ ||\xi||G\psi_{\phi} < \infty. \ $ \par

\vspace{4mm}

 \  {\bf  Our aim in this short  report is description of the associate and dual spaces for the Grand Lebesgue Spaces.} \par

\vspace{4mm}

 \ We intend to simplify one in the articles [6]-[7], [17]. \par

\vspace{4mm}

\section{ Main result: structure of associate spaces. }

 \vspace{4mm}

 \ Recall first of all that the {\it associate} space $ \  Y' \ $ for arbitrary rearrangement one  $ \  Y \ $ builded over source measurable space
$ \  (X,B, \mu)  \  $ consists by definition on all the {\it linear} functionals of the forms

$$
l(f) = l_g(f) \stackrel{def}{=} \int_X f(x) \ g(x) \ \mu(dx), \ f(\cdot) \in Y, \eqno(2.1)
$$
for (certain) measurable function $ \ g: X \ \to R. \ $ It will be presumed that this functional  $ \ l_g(\cdot) \ $ is bounded:

$$
||l_g||Y' \stackrel{def}{=} \sup_{f: ||f||Y \le 1} |l_g(f)| < \infty. \eqno(2.2)
$$
 \ The detail description of these spaces is explained in the classical monograph of C.Bennet and R.Sharpley [3], chapters 1,2. \par

 \ We will identify in this case as usually for convenience the linear functional $ \ l_g \ $ with its generating function $ \ g: \ $

$$
||g||Y' := ||l_g||Y'.
$$

\vspace{4mm}

 \ Let now $ \ Y = G\psi \ $ for some $ \ \psi = \psi(\cdot) \in \Psi(a,b), \ $ where as above  $ \ a = \const \ge 1, \ b = \const \in (a, \infty]. \ $
 Let also $ \ 0 \ne f \in G\psi; \ $ one can suppose without loss of generality $ \ ||f||G\psi= 1. \ $  Therefore

$$
|f|_p \le \psi(p), \ p \in (a,b).
$$

 \  Denote as usually

 $$
 q= q(p) = p' := p/(p-1), \ a' = a/(a-1), \ b' = b/(b-1),
 $$
so that $ \ \infty' = 1. \ $  Introduce the so-called  {\it adjacent} function $ \ \nu = \nu(q) = \nu[\psi](q), \ q \in (b', a') \ $ as follows

$$
\nu(q) = \nu[\psi](q) \stackrel{def}{=} \frac{1}{\psi(q/(q-1))}. \eqno(2.3)
$$

 \ Let us estimate the functional $ \ l_g(f) \ $ from the relation (2.1). We apply the classical H\"older's inequality for the values correspondingly $ \ p  \in (a,b) \ $
 and $ \ q = q(p) \in (b',a') $

$$
|l_g(f)| \le |f|_p \ |g|_{q(p)}  \le \psi(p) \ |g|_q = \frac{|g|_q}{\nu[\psi](q)}.
$$

\ Define a following functional

$$
V(g) = V[\psi](g) := \inf_{q \in (b', a')} \left[ \ \frac{|g|_q}{\nu[\psi](q)} \  \right], \eqno(2.4)
$$

 \ We  proved actually the following proposition:

\vspace{4mm}

 \ {\bf Theorem 2.1. } Assume that $ \  V(g) \ < \infty; \ $ then $ \ g \in G\psi' \ $ and wherein  $ \ ||g||G\psi' \le V(g). \ $ \par
 \  As a  slight consequence:

$$
  G\nu \subset G\psi'; \ ||l_g||(G\psi)' \le ||g||G\nu[\psi]. \eqno(2.4a)
$$

\vspace{4mm}

 \ {\bf Example 2.1.} Let the measure $ \ \mu \ $ be probabilistic: $ \ \mu(X) = 1. \ $ If $ \ \psi(p) := \psi_{(r)}(p) = 1, \ p \in [1,r], \ $ where
 $ \ r = \const \in [1, \infty), \ $ so that $ \ a = 1, \ b = r, \ \supp \ \psi_{(r)} = [1,r], \ $
 then the space $ \ G\psi_{(r)} \ $ coincides  with the classical Lebesgue - Riesz space $ \  L_r = L_r(X).  \  $ \par
  \ In  this case

 $$
 V(g) = |g|_{r'}, \ r' = r/(r-1).
 $$
 \ Thus, in this example the associate space is quite equal to the its conjugate. \par

\vspace{4mm}

\ {\bf Example 2.2.} Let now again in the probabilistic case, i.e. when $ \ \mu(X) = 1 \ $

$$
\psi(p) = \psi_m(p) \stackrel{def}{=} p^{1/m}, \ p \in [1,\infty), \ m = \const > 0.
$$
 \ The case $ \ m = 2 \ $ correspondent to the subgaussian case. \par
 \ We have

$$
\nu_m(q) \stackrel{def}{=} \nu \left[\psi_m \right](q)   = \left( \frac{q-1}{q} \right)^{1/m};
$$
therefore

$$
||g||(G\psi_m)' \le \inf_{q \in (b',a')} \left\{ \ \left( \ \frac{q}{q-1} \ \right)^{1/m} \ |g|_q \  \right\}.
$$

 Of course

$$
\nu_m(q) \asymp (q-1)^{1/m}, \ q \in (1,2); \ \nu_m(q) \asymp  1, \ q \in [2, \infty).
$$
 \ Notice that in the considered case the  the estimate (2.4a) gives a weak result. Namely, as long as for any r.v. $ \ \eta \ $

$$
 \sup_{q \ge 2} |\eta|_q = \lim_{q \to \infty} |\eta|_q = \vraisup_{x \in X} |\eta(x)| = |\eta|_{\infty},
$$
we deduce on the basis of relation (2.4a) a trivial estimate

$$
||g||(G\psi_m)' \le C(m) |g|_{\infty}, \ C(m) \in (0,\infty).
$$

\vspace{4mm}

\section{ Main result: structure of dual spaces. }

\vspace{4mm}

 \ We need to introduce the following addition notations, conditions and definitions.

 $$
h(p) =  h[\psi](p) := p \ \ln\psi(p),  \ p \in (a,b),
 $$

$$
 V(u) = V[\psi](u):=  h^*[\psi](\ln |u|), \ \psi \in \Psi.
$$
  Define also the following Young-Orlicz function

$$
\  N[\psi](u) := \exp( V(u) ) =   \exp( V[\psi](u) )  :=
$$

$$
 \exp \left[ \ h^*[\psi](\ln |u|) \ \right], \ u \ge e; \ N[\psi](u) = C \ u^2, \ |u| < e. \eqno(3.1)
$$

 \ Recall that the well - known Young-Fenchel, or Legendre transform $ \ h^* \ $  for the real valued function $ \ h = h(z),  \ $ not necessary to be convex,
  is defined as follows

$$
h^*(v) \stackrel{def}{=} \sup_{z \in \Dom(h) } ( v z - h(z)).
$$

 \ Further,  denote by $ \ \Gamma = \{ \ \gamma  \ \}  \ $ the collection of all {\it  finite-additive } numerical valued
functions defined on the sigma - field $ \ B. \ $ \par

 \ Of course, the finite additivity does  not exclude the "ordinary" countably sigma - additivity. \par

 \ Define following [18], see also [19],  norm on the set $ \ \Gamma: \ $

$$
|||\gamma|||  = |||\gamma|||_{\psi} \stackrel{def}{=} \sup \left\{ \ \int_X f(x) \ \gamma(dx): \ ||f||G\psi \le 1 \ \right\}, \eqno(3.2)
$$
where the integral in (3.2) may be understood in the sense of R.G.Bartle [2]; see also [18]. On the other hands, one can  choose only
the function $ \ f(\cdot) \ $ to be {\it simple}, or on the other words, {\it step function}, i.e. taking only finite possible values on the sets of finite "measure" $ \ \gamma: \ $ \par

$$
f(x) = \sum_{i=1}^n c_i \ \chi_{D(i)}(x), \ D(i) \in B, \ c_i \in R,  \ n  = 1,2, \ldots, \ \gamma(D(i)) \in R;
$$
where as ordinary $ \ \chi_D(x) \ $ denotes an indicator function of the measurable set $ \ D: \ \chi_D(x) = 1, \ x \in D; \ \chi_D(x) = 0, \ x \notin D. \ $ \par

 \ Obviously, for such a functions

$$
\int_X f(x) \ \gamma(dx) = \sum_{i=1}^n c_i \ \gamma(D(i)).
$$

\vspace{4mm}

{\bf Theorem 3.1}   Assume  in addition to the foregoing conditions imposed on the function $  \ \psi(\cdot) \ $
 that the function $  \  V(x) = V[\psi](x)  \  $ satisfies the following restriction:
$$
\exists \ \alpha = \const \in (0,1), \ \exists K  = \const > 1, \forall  x \in (0,\infty) \ \Rightarrow
V(x/K) \le \alpha \ V(x), \eqno(3.3)
$$
which is in turn a slight analog of the famous $ \ \Delta_2 \ $ condition.\par

 \ We assert that arbitrary continuous linear functional $ \ l(f), \ f \in G\psi \ $ has an unique representation of the form

$$
l(f)  = l^{\gamma}(f) = \int_X \ f(x) \ \gamma(dx), \ \gamma \in \Gamma, \ |||\gamma|||_{\psi} < \infty, \eqno(3.4)
$$
and conversely each functional $ \  l^{\gamma}(\cdot)  \  $ of the form (3.4) belongs to the conjugate (dual) space $ \ (G\psi)^* \ $ and moreover

$$
|| l^{\gamma}||(G\psi)^* = |||\gamma|||_{\psi}; \eqno(3.5)
$$
on the other words, the (Banach) space $ \  S = S(\psi) :=(\Gamma, \ ||| \ \cdot \ ||| ) \ $ coincides with  the dual (conjugate) space $ \  (G\psi)^* \ $ relative the equality (3.5). \par

 \vspace{4mm}

 \  {\bf Proof.} {\bf I.}  Suppose $ \  \gamma \in S  \  $ and following $ \ |||\gamma|||_{\psi} < \infty. \ $ The equality (3.5) follows immediately from the direct
 definition of the spaces $ \ S(\psi); \ $  obviously, it is sufficient to verify (3.5) for all the step functions. \par
  \ It is worth to note that this conclusion is true still without the condition  (3.3) \par

 \vspace{4mm}

  \ {\bf 2.} Conversely, let the space  $  \  G[\psi], \ \psi \in \Psi(a,b),  \ 1 \le a < b \le \infty  \  $  be given. Assume also that the condition (3.3)  is satisfied. \par
 \ It is proved in particular in  [12] that the Grand Lebesgue Space $ \  G\psi \ $ in this case quite coincides with the so-called {\it exponential} Orlicz space $ \  L(N[\psi]) \ $
 builded over our measurable triplet $ \ (X,B,\mu) \ $ equipped with the correspondent Young-Orlicz function

$$
N[\psi](u) :=  \exp( V(u) ) =   \exp(V[\psi](u)).
$$
 \ The norm in this space may be used the ordinary Luxemburg norm. \par
 \ The conjugate spaces for these ones are calculated in [18]; see also [19]. They have the form  (3.5). \par
 \ This completes the proof of theorem 3.1. \par

\vspace{4mm}

 \ {\bf Remark 3.1.} The condition (3.3) is satisfied for very popular class of  $ \ G \psi \ $ spaces, indeed, if for instance

$$
V(x) = C \ x^m, \ C,m = \const \in (0,\infty), \ x > 0,
$$
as well as in the case when

$$
V(x) = C \ x^m \ L(x), C,m = \const \in (0,\infty), \ x > 0,
$$
where $ \ L = L(x) \ $ is positive continuous slowly varying simultaneously $ \ x \to 0+ \ $ and as $ \ x \to \infty \ $ function
such that

$$
L(x/K) \le \alpha \ K^m \ L(x), \ x > 0. \eqno(3.6)
$$

 \ In turn, the last condition is satisfied if the function $ \ L(x) \ $ is in addition strictly  increasing. \par

\vspace{4mm}

\section{ Again about associate space.}

\vspace{4mm}

 \ Let us return to the problem of finding of associate space to the Grand Lebesgue ones. We intend to apply the described above approach
throughout embedding into Orlicz spaces  [13], [14], [17]. \par
 \ In detail, let $ \ \mu(X) = 1 \ $ and let the function $ \ \psi(\cdot) \ $ satisfy the condition (3.3). Suppose
 also that the function $ \ f:  X \to R \  $ be from the space $ \ G\psi; $ then it belongs also to the correspondent Orlicz's space $ \ L(N[\psi])   \ $ and herewith

$$
||f||L(N[\psi]) \le C(\psi)  \ ||f||G\psi, \  C(\psi) < \infty,
$$
and inverse inequality is also true.\par
 \ One can use the famous H\"older's inequality

$$
|l_g(f)| = \left| \ \int_X f(x) \ g(x) \ \mu(dx) \  \right| \le  2 \ ||f||L(N[\psi]) \ ||g||L(N^*[\psi]) \le
$$

$$
2 \ C(\psi) \ ||f||G\psi \ ||g||L(N^*[\psi]),
$$
see e.g. [9], [19], [20]. \par

\vspace{4mm}

  \ To summarize: \par

 \vspace{4mm}

 \ {\bf Theorem 4.1.}  Suppose once again that the function $ \ \psi(\cdot) \ $ is from the set $ \ \Psi \ $
 and satisfies the condition (3.3).  The Orlicz space $ \ L(N^*[\psi])  \  $ over our measurable space $ \ (X, B, \ \mu) \ $  quite coincides
to the associate  $ \ (G\psi)' \ $ for the Grand Lebesgue Space $ \ G\psi \ $ and moreover

$$
||l_g||(G\psi)' \le 2 \ C(\psi) \ ||g||L(N^*[\psi]). \eqno(4.1)
$$

\vspace{4mm}

 \ {\bf An example.} \par

 \ Let $ \ \mu(X) = 1; \ \psi(p) = \psi_2(p) := p^{1/2}, \  - \ $ a subgaussian case, see e.g. [5].  Then the correspondent conjugate Young-Orlicz function has a form

$$
N^*[\psi_2](y) \asymp  |y| \ \ln^{1/2}(e + |y|), \ |y| \ge 1.
$$
 \ A slight  generalization: $ \ \psi(p) = \psi_m(p) := p^{1/m}, \ m = \const > 0. \ $ Then

$$
N^*[\psi_m](y) \asymp  |y| \ \ln^{1/m}(e + |y|), \ |y| \ge 1.
$$

 \ Many other examples of calculation of these functions $ \  N^*[\psi](p), \  $ including the cases when $ \ \psi(p) = p^{1/m} \ L(p), \ p \in [1, \infty), \ $
 where $ \ L = L(p) \ $ is positive continuous slowly varying as $ \ p \to \infty \ $ function; $ \ \psi(p) = \exp \left(C p^{\beta} \right),  \ C, \beta, p = \const \in (0,\infty) \ $ etc,
 may be found in the articles [12], [14]. \par

\vspace{6mm}

\section{Concluding remarks.}

\vspace{4mm}

 {\bf A. } It is no hard by our opinion to generalize obtained here results into a multidimensional case of  $ \  G\psi \ $ spaces.\par

\vspace{4mm}

 \ {\bf B. } As long as the $ \ B(\phi) \ $ spaces are particular cases of Grand Lebesgue ones, we obtained on the way also the associate and dual for these
spaces, as well. \par

\vspace{4mm}

\ {\bf C. } It is interest by our opinion to extend the integral representation for dual and associate linear continuous functionals for many wide classes of
rearrangement invariant spaces.\par

 \begin{center}

 \vspace{6mm}

 {\bf References.}

 \vspace{4mm}

\end{center}

 {\bf 1. Rodrigo Banuelos and Adam Osekovski.} {\it  Weighted square function estimates.} \\
  arXiv:1711.08754v1  [math.PR]  23 Nov 2017 \\

\vspace{3mm}

{\bf 2.  R.  G.  Bartle.} {\it  A   general  bilinear   vector  integral.}
   Studia  Math.  15   (1956),  337-352.\\

\vspace{3mm}

{\bf 3. Bennet C., Sharpley R.}  {\it  Interpolation of operators.} Orlando, Academic
Press Inc., (1988). \\

 \vspace{3mm}

{\bf 4.  Buldygin V.V., Kozachenko Yu.V. }  {\it Metric Characterization of Random
Variables and Random Processes.} 1998, Translations of Mathematics Monograph, AMS, v.188. \\

 \vspace{3mm}

{\bf 5. X.Fernique.} {\it Regularite des trajectoires des fonctions aleatoires gaussiennes.} Lectures Notes in Mathematics, {\bf 480}, (1975).
Springer Verlag, Berlin-Heidelberg.\\

\vspace{4mm}

 {\bf 6. A. Fiorenza.}   {\it Duality and reflexivity in grand Lebesgue spaces. } Collect. Math.
{\bf 51,}  (2000), 131-148. \\

 \vspace{3mm}

{\bf 7. A. Fiorenza and G.E. Karadzhov.} {\it Grand and small Lebesgue spaces and
their analogs.} Consiglio Nationale Delle Ricerche, Instituto per le Applicazioni
del Calcoto Mauro Picone”, Sezione di Napoli, Rapporto tecnico 272/03, (2005).\\

 \vspace{3mm}

 {\bf 8. Dudley R.M.} {\it Uniform Central Limit Theorem.} Cambridge, University Press, (1999), 352-367.\\

\vspace{3mm}

 {\bf 9. Daniele Imparato.} {\it  Martingale inequalities in exponential Orlicz spaces.}
 Journal of inequalities in pure and applied mathematic.
 Volume 10  (2009), Issue 1, Article 1,  pp. 1-10.\\

\vspace{3mm}

{\bf 10.  T. Iwaniec and C. Sbordone.} {\it On the integrability of the Jacobian under minimal
hypotheses. } Arch. Rat.Mech. Anal., 119, (1992), 129-143. \\

 \vspace{3mm}

{\bf 11. Kozachenko Yu. V., Ostrovsky E.I. }  (1985). {\it The Banach Spaces of random Variables of subgaussian Type. } Theory of Probab.
and Math. Stat. (in Russian). Kiev, KSU, 32, 43-57. \\

\vspace{3mm}

{\bf 12. Kozachenko Yu.V., Ostrovsky E., Sirota L.}  {\it Relations between exponential tails, moments and
moment generating functions for random variables and vectors.} \\
arXiv:1701.01901v1 [math.FA] 8 Jan 2017 \\

\vspace{3mm}

{\bf 13. Ostrovsky E.I. } (1999). {\it Exponential estimations for Random Fields and its
applications,} (in Russian). Moscow-Obninsk, OINPE. \\

 \vspace{3mm}

{\bf 14. Ostrovsky E. and Sirota L.} {\it Vector rearrangement invariant Banach spaces
of random variables with exponential decreasing tails of distributions.} \\
 arXiv:1510.04182v1 [math.PR] 14 Oct 2015 \\

 \vspace{3mm}

{\bf 15. Ostrovsky E.I.} {\it Generalization of norm of Buldygin-Kozachenko and Central Limit Theorem in Banach space. }
Probability Theory Applications, 1982, V. 27 , Issue 3, p. 618. \\

 \vspace{3mm}

{\bf 16. Ostrovsky E., Rogover E. } {\it Exact exponential bounds for the random field
maximum distribution via the majorizing measures (generic chaining).} \\
 arXiv:0802.0349v1 [math.PR] 4 Feb 2008 \\

\vspace{3mm}

 {\bf 17. Ostrovsky E. and Sirota L. }  {\it  Bilateral small Lebesgue spaces.} \\
arXiv:0904.4535v1  [math.FA]  29 Apr 2009 \\

 \vspace{3mm}

{\bf 18.  M.M.Rao.} {\it Linear functional on Orlicz spaces.} Pacific Journal of Mathematics,
Vol.  25,  No. 3,  535-585,  1968. \\

\vspace{3mm}

{\bf 19. M.M.Rao, R.Z.Ren.} {\it  Theory of Orlicz spaces.} Marcel Dekker. \ Berlin, Heidelberg, 1992. \\

\vspace{3mm}

{\bf 20.  O.I.Vasalik, Yu.V.Kozachenko, R.E.Yamnenko.} {\it $ \ \phi \ - $ subgaussian  random processes. } Monograph, Kiev, KSU,
2008;  (in Ukrainian). \\

\vspace{3mm}

\end{document}